\pdfoutput=1 
\documentclass[10pt]{article}

\usepackage{amsmath,amssymb,amsfonts}
\usepackage{amsthm}
\usepackage{hyperref}
\usepackage[T1]{fontenc}
\usepackage{lmodern}
\usepackage{listings}
\usepackage{upquote} 
\usepackage{cmap}    
\usepackage{xcolor}
\usepackage{tikz}
\usepackage{caption}
\usepackage{float}
\usepackage{graphicx}

\lstdefinelanguage{Magma}{
  morekeywords={intrinsic,function,return,for,in,to,do,end,if,then,else,elif,while,repeat,until,require,assert,is,where,declare,assign,print,procedure,true,false,seqenum,SeqEnum,Prj,PrjProd,Aff,RngIntElt,RngMPolElt,Sch,Pt,ModRng},
  sensitive=true,
  morecomment=[l]{//},
  morecomment=[s]{/*}{*/},
  morestring=[b]"
}
\title{The Linear System Package of Magma}
\author{Carlos Rito}
\date{}

\begin{document}
\maketitle


\begin{abstract}
We present a complete reimplementation of the LinearSystem package of Magma,
with substantial improvements in design and performance.
The resulting efficiency enables computations that were previously out of reach.
We briefly describe the design principles, capabilities, and algorithms of the new implementation
and illustrate them with examples that showcase its power.
Rather than comparing speeds, our goal is to advertise
the package by demonstrating what can now be achieved in practice.
We also add one core capability: computing linear systems of plane curves with prescribed non-ordinary singularities.


\end{abstract}

\section{Introduction}

A \emph{complete linear system} on a variety is the projective space consisting of all
effective divisors that are linearly equivalent to a given divisor.
A \emph{linear system} is any projective linear subspace of this space.

Linear systems are fundamental across Algebraic Geometry: they provide embeddings and birational maps, govern pluricanonical models, allow the study of subvarieties.

Magma \cite{BCP} has long included a package for computations with linear systems,
originally written by Gavin Brown and Paulette Lieby about 25 years ago.
It offers a rich and broadly useful feature set. In this paper we introduce a new implementation---written from scratch---whose design goal is speed while retaining the standard functionality. In practice, essentially all operations are much faster; in particular, linear systems through points are routinely more than 500 times faster. Tasks that previously tended not to finish (e.g. imposing many thousands of conditions) now become routine.

We do not attempt to analyze the causes of these speed differences. Instead, we outline the guiding principles of the new implementation and emphasize what it enables. The central design choice is to keep objects light at creation time and to materialize heavy data only on demand.

Beyond performance, the package adds a new capability: the computation of plane curves with \emph{non-ordinary} singularities (such as cusps, tacnodes, and higher contact). These conditions are enforced algorithmically by tracking tangent directions along blowups and translating them into linear constraints, so the entire process remains linear-algebraic and scales well.

The goal of this paper is to make the community aware of what can be done with these tools. We give a brief account of the internal data model and constructors, explain how to impose geometric conditions (through points or subschemes), and then focus on examples that we hope the reader will find both fun and powerful, including finite-field constructions of quintic surfaces with many nodes or cusps.

Section~2 presents the data model and constructors (complete systems, systems from sections, and from matrices/monomials), together with on-demand coefficient maps. Section~3 covers restrictions and imposed conditions (ordinary multiplicities at points, containment of subschemes in projective vs.\ affine ambients, and fast workflows for images, parameter loci, and parameter recovery). Section~4 introduces the new machinery for non-ordinary plane singularities via blowups and tangent directions. Section~5 showcases applications to singular quintic surfaces, including \(\mathbb{Z}/5\)– and \(\mathbb{Z}/6\)–invariant families yielding many nodes or cusps.

Some of the new reimplementations entered Magma in version V2.28-1; the complete reimplementation of the package is available from version V2.28-16 onward.

Readers without a Magma license can still run most of the Magma code in this paper via the Magma Online Calculator \cite{MagmaCalc}.

\subsubsection*{Acknowledgments}

The author thanks John Cannon and John Voight, and the Sydney Mathematical Research Institute (SMRI), for the invitations to visit Sydney.

The author warmly thanks Allan Steel for promptly improving Magma’s speed in several important situations.

This research was partially financed by Portuguese Funds through FCT
(Funda\c c\~ao para a Ci\^encia e a Tecnologia) within the Project UID/00013:
Centro de Matem\'atica da Universidade do Minho (CMAT/UM).

\section{Data model and constructors}

This section describes the internal representation of a linear system in Magma,
together with the basic constructors available to the user. 
The guiding principle is: \emph{keep the object as light as possible at the time of definition, 
and compute heavy data (bases, matrices, maps) only when needed.}

\subsection{The \texttt{LinearSys} object}

A linear system is an object of type \texttt{LinearSys} whose key attributes include:
\begin{itemize}
  \item \texttt{Ambient} (ambient space, projective/affine/product),
  \item \texttt{Degree} (an integer degree, or a multidegree sequence for multigraded ambients),
  \item \texttt{Sections} (optional, an explicit list of polynomials),
  \item \texttt{Monomials} and \texttt{Matrix} (optional, a coefficient matrix together with the monomial list),
  \item \texttt{IsComplete}, \texttt{Echelonized}, \texttt{IndependentSections} (bookkeeping flags),
  \item \texttt{CoefficientSpace}, \texttt{CoefficientMap} and \texttt{PolynomialMap} (created on demand).
\end{itemize}

There are two parameters:
\begin{itemize}
  \item \texttt{CheckBasis} (default \texttt{true}): whether to verify linear independence
        and, if necessary, switch to the matrix form by echelonizing the coefficient matrix;
  \item \texttt{ChangeBasis} (default \texttt{false}): when \texttt{true}, even linearly independent
        inputs are replaced by an echelonized basis.
\end{itemize}

\subsection{Constructors}

The possibilities are:

\medskip
\noindent
\texttt{LinearSystem(Ambient,Degree)}
\newline
\texttt{LinearSystem(Ambient,Sections)}
\newline
\texttt{LinearSystem(Ambient, Matrix, Monomials)}

\medskip

A complete system on a projective space of degree $d$ is initially stored
as the ambient plus the integer $d$ (and a flag \texttt{IsComplete}).
A basis of sections or a coefficient matrix is only materialized when required
(e.g.\ querying \texttt{Sections}, \texttt{Dimension}, applying restrictions, etc.).

\medskip

For example, from matrix and monomials:
\begin{lstlisting}[language=Magma]
P2<x,y,z>:=ProjectiveSpace(Rationals(),2) ;
mon:=[x^2,y^2,z^2,x*y,x*z];
M:=Matrix([
[1,0,1,0,0],
[0,1,0,0,-1],
[0,1,0,1,0],
[0,0,0,0,1]
]);
L:=LinearSystem(P2,M,mon);
\end{lstlisting}

Sections are created (and then stored) on demand:

\begin{lstlisting}[language=Magma]
Sections(L) eq [x^2+z^2,y^2-x*z,x*y+y^2,x*z];
\end{lstlisting}

A linear system can be created from the sections:

\begin{lstlisting}[language=Magma]
J:=LinearSystem(P2,Sections(L):ChangeBasis:=true);
Sections(J) eq [x^2+z^2,x*z,y^2,x*y];
\end{lstlisting}
But note that in this case the sections have been echelonized, so we are considering a linear system with a different basis.

\medskip

With \texttt{CheckBasis:=true} the constructor builds a coefficient map against a monomial list,
echelonizes and stores the matrix+monomial form.
If one already knows the list is a basis and wish to avoid preprocessing, set
\texttt{CheckBasis:=false}; then the sections are stored verbatim and linear algebra
is deferred until a computation requires it. Setting \texttt{ChangeBasis:=true}
forces an echelonized form even for independent inputs.

Echelonization may reduce the number of monomials, but it can introduce larger coefficients when working over the rationals. Whether to use it depends on the specific context.

Those two options were not available in the previous implementation of the package.

\medskip

Variants exist for affine ambients and for products (where multidegree is enforced).

\subsection{Coefficient maps}

Two natural maps are associated to every system:
\medskip

The {\bf coefficient map}, sending a section to its coordinate vector in the chosen basis.
The previous implementation always computed a coefficient map at creation,
ensuring fast later use but making creation itself very slow.
In the new version, the map is built only on demand,
so creation is immediate while the cost is deferred to when/if it is actually needed;
\medskip

The {\bf polynomial map}, the inverse map from coefficient vectors
        to polynomials.
\medskip

These maps allow a transfer between the geometric and linear algebraic viewpoints.

\medskip

For {\bf example}, suppose we have a large sequence of polynomials $s$ and need an efficient way to compute the coefficients of given polynomials in terms of $s$ --- for instance, when this must be done thousands of times. We define the linear system $L$ given by $s$, keeping the basis unchanged. Even if $s$ is not linearly independent, this still works (though the coefficients are not unique). The computation of the coefficient map of $L$ may be time-consuming, but once it is available, applying it is very fast.
\medskip

 Let's consider a sequence of 100 polynomials of degree 50 in 4 variables:

\begin{lstlisting}[language=Magma]
P3:=ProjectiveSpace(Rationals(),3);
L:=LinearSystem(P3,50);
s:=[Random(L,[-10..10]):i in [1..100]];
Ls:=LinearSystem(P3,s:CheckBasis:=false);
f:=Random(Ls,[-10..10]);
h:=CoefficientMap(Ls);  // 26 sec
cfs:=h(f);              // 0.2 sec

\end{lstlisting}

\medskip

We note that membership testing (\texttt{f in L}) computes (and stores) the coefficient map, if not computed before.

\subsection{Reduction, base scheme and trace}

So far we have considered linear systems defined on ambient spaces.
In practice the relevant situation is to restrict these systems
to a given variety $X \subseteq \mathbb{A}^n$ or $\mathbb{P}^n$.
 
A linear system $L$ is given by a family of sections,
and it is important to understand the common zero locus of these sections.
This is obtained in Magma using \texttt{BaseScheme(L)}.
It consists of the points of the ambient space
where all members of $L$ vanish simultaneously.

If all sections of $L$ share a fixed component, it can be removed using\newline
\texttt{Reduction(L)}, leaving only the moving part of the system.

To restrict a linear system $L$ on the ambient space to a variety $X$,
one discards all sections that vanish identically on $X$.
This is accomplished with the command \texttt{LinearSystemTrace(L,X)}.
The result is the linear system induced on $X$, containing precisely
the sections of $L$ that cut non-trivial divisors on $X$.
This is equal to \texttt{Complement(L,LinearSystem(L,X))}.

\medskip For {\bf example}, if we take a random surface $S$ cut out by 4 quadrics
in $\mathbb P^6,$ its bicanonical system is given by all quadrics that do not vanish identically
on $S$. Its dimension is 23.
\begin{lstlisting}[language=Magma]
P6:=ProjectiveSpace(Rationals(),6);
L2:=LinearSystem(P6,2);
S:=Surface(P6,[Random(L2,[1..10]):i in [1..4]]);
T:=LinearSystemTrace(L2,S);
Nsections(T) eq 24;
\end{lstlisting}

\medskip We note that \texttt{Nsections(T)} and \texttt{\#Sections(T)} return the same value,
but the former avoids computing the full list \texttt{Sections(T)}.

\section{Restrictions and imposed conditions}\label{sec:restrictions}

A large part of practical work with linear systems consists of imposing
geometric conditions: passing through points (with multiplicities), containing a given
subscheme, or enforcing prescribed singularities. This section describes the
interfaces and the underlying algorithms for these tasks.

\subsection{Through points with ordinary multiplicities}

Let $L$ be a linear system on an ambient $A$ (projective or affine). Given points
$p_1,\dots,p_r$ and nonnegative integers $m_1,\dots,m_r$, we can compute
the subsystem of members whose multiplicity at $p_i$ is at least $m_i$ for each $i$.

Write the sections of $L$ as $s_1,\ldots,s_N$.
The condition “$F\in L$ has multiplicity $\ge m$ at $p$” is linear
in the coefficients of $F=\sum a_js_j$ and is enforced by the vanishing of all
partial derivatives of order $< m$ at $p$. The implementation builds an
evaluation/derivative matrix and
extracts a basis of the nullspace. This uses only linear algebra,
thus it is fast, allowing the computation of systems through thousands of points.

\medskip
For {\bf example}, over a finite field and through 3275 points:

\begin{lstlisting}[language=Magma]
K:=GF(397); P:=ProjectiveSpace(K,3);
pts:=[P![Random(K):i in [1..4]] : j in [1..3275]];
L:=LinearSystem(P,25);
J:=LinearSystem(L,pts);   // 1.2 sec
Nsections(J) eq 1;
\end{lstlisting}

\medskip
We compute a plane curve of degree 20 with many ordinary singularities, over the rationals:

\begin{lstlisting}[language=Magma]
A:=AffineSpace(Rationals(),2);
m:=[2,2,2,2,2,2,3,3,3,3,3,5,5,5,7,7,8,9];
pts:=[A![Random(1,40),Random(1,40)]:i in [1..#m]];
L:=LinearSystem(A,20);
J:=LinearSystem(L,pts,m);   // 0.6 sec
Nsections(J) eq 1;
\end{lstlisting}

\subsection{Containing a subscheme}\label{subsec:containX}

Let $X\subset A$ be a subscheme of the ambient. The subsystem of $L$ whose members contain $X$
is computed by \texttt{LinearSystem(L,X)}. The approach depends on whether
the ambient is projective or affine.

\medskip
\noindent\textbf{Projective ambient.}
To avoid problems with the irrelevant ideal, we start by saturating the ideal $I$ of $X.$\newline
(For example if the ideal of $X\subset\mathbb P^1$ is generated by $(x^2,xy),$
then the saturation of $I$ is generated by $(x)$).

Then for each generator $q$ of $I,$ we form all products $q m$ with monomials $m$ so that
$\deg(q m)=\deg(L)$ and collect them into a candidate list of sections.
We then perform a coefficient-space elimination (via echelonization and nullspaces) to
extract a basis of the subspace of $L$ vanishing on $X$.
Only linear algebra involved in these computations.

\medskip
For {\bf example}, we take a surface $X$ in the 5 dimensional projective
space and compute all polynomials of degree 15 that vanish on $X.$
One may choose to mark the ideal of $X$ as already saturated,
so that saturation is not recomputed (saturation can be computationally expensive).

\begin{lstlisting}[language=Magma]
K:=Rationals();
P:=ProjectiveSpace(K,5);
L:=LinearSystem(P,15);
s:=[Random(L,[1..10]):i in [1..3]];
X:=Scheme(P,s:Saturated:=true);
J:=LinearSystem(L,X);
\end{lstlisting}

\medskip
\noindent\textbf{Affine ambient.}
Here the previous approach doesn't work because monomials of the same degree may cancel and give rise
to polynomials of lower degree. The solution to the {\em ideal membership problem} is solved via Gr\"obner bases.

Let $s_1,\dots,s_n$ be the sections of $L$ and write an
unknown $F=\sum_1^n a_i s_i$ where the $a_i$ are new coefficient variables.
We extend the polynomial ring of $X$ by adding the variables $a_i$
and compute the Gr\"obner basis $G$ of $X$ in this new ideal.
After computing the normal form \texttt{NF:=NormalForm(F,G)},
we just need to consider the $a_i$'s such that $NF = 0$.

\subsection{Some practical applications}

In practice, three scenarios are particularly effective.

\medskip
\noindent\textbf{Images.}
Let $X$ be a variety and $f\colon X\dashrightarrow \mathbb{P}^n$ a rational map. Computing the ideal of $f(X)$ via Gr\"obner-basis elimination is often prohibitively expensive. For a fixed degree $d$, the command \texttt{ImageSystem(f,X,d)} returns the linear system of degree $d$ hypersurfaces containing $f(X)$. Because it uses only linear algebra, this is typically faster than Gr\"obner bases. However, for large instances it can still be a bottleneck. A faster alternative is to proceed as follows:

\begin{description}
  \item[$\cdot$] Working over finite fields $\mathbb{F}_p$, sample many points on $X$ (e.g.\ by intersecting with random hyperplanes) and evaluate $f$ to obtain many points of $f(X)$.
  \item[$\cdot$] Define a linear system $L$ of degree $d$ and impose the point conditions with \texttt{LinearSystem(L,pts)} to recover all degree $d$ polynomials vanishing on $f(X)$.
  \item[$\cdot$] Increase $d$ and aggregate equations until stabilization is observed.
  \item[$\cdot$] Optionally, repeat across several primes and lift the resulting coefficients to characteristic $0$ (e.g. via CRT/RationalReconstruction in Magma).
\end{description}

Since one can impose thousands of point conditions quickly, this approach is highly effective for describing $f(X)$ by equations while avoiding costly elimination.

\medskip
\noindent\textbf{Families -- parameter space.}
Suppose a construction associates to each point $u\in S\subset\mathbb{P}^n$ a variety $X_u$,
and we can \emph{produce many parameter points} $u$ for which $X_u$ exists (or meets prescribed properties).
Here $S$ is the unknown parameter locus we wish to recover.
As above, we sample many such parameters $u_1,\dots,u_m\in S$ (typically over finite fields)
and compute \texttt{LinearSystem(L,pts)} to recover all degree-$d$ polynomials vanishing at the $u_i$.
Repeating for other values of $d$, we eventually obtain the full set of defining equations of $S$.

\medskip
\noindent\textbf{Parameter recovering.}
If, for sampled parameters $u\in\mathbb{P}^n$, we can compute the \emph{defining equations} of $X_u$, we can recover the \emph{parameter–dependence} of their coefficients.
Say that one of these equations is
\[
   x^d + \frac{N_1(p)}{D_1(p)}x^{d-1}y + \cdots + \frac{N_m(p)}{D_m(p)}w^d.
\]
The challenge is to recover each rational function $\frac{N_i(p)}{D_i(p)}$. 
We introduce a new variable $v$ and compute many points
        $(\overline{p}, \overline{v})$ such that
        $\overline{v}=\tfrac{N_i(\overline{p})}{D_i(\overline{p})}$.
Then use linear systems through these many points to recover the polynomial
$F(p,v) = N_i(p) - D_i(p) v.$ This gives the desired rational function.

\section{Plane curves with non-ordinary singularities}\label{sec:plane}

A new feature of the package is the ability to impose \emph{non-ordinary} singularities
on affine plane curves. 
Given a plane curve defined by $F(x,y)=0$,
an ordinary singularity at $p$ is enforced by requiring that
$F$ and its partial derivatives up to order $m-1$ vanish at $p$.
Non-ordinary singularities (cusps, tacnodes, higher contacts) require
tracking tangent directions through blowups. 
For example, a tacnode is resolved after blowing up once and then
requiring multiplicity two at the infinitely near point determined by the tangent direction.

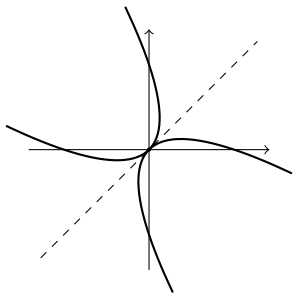
\begin{figure}[h]
\centering
\begin{tikzpicture}[scale=0.8]
  \draw[->] (-2,0)--(2,0) node[right] {};
  \draw[->] (0,-2)--(0,2) node[above] {};

  \draw[dashed] (-1.8,-1.8) -- (1.8,1.8) node[above right] {};

  \begin{scope}[rotate=45]
    \draw[domain=-1.4:1.4,smooth,thick,variable=\x] plot ({\x},{\x*\x});
    \draw[domain=-1.4:1.4,smooth,thick,variable=\x] plot ({\x},{-\x*\x});
  \end{scope}

  \fill (0,0) circle (1pt);
  \node[below right] at (0,0) {};
\end{tikzpicture}
\caption*{Tacnode with tangent direction $(1,1)$.}
\end{figure}

Let $S_1 \to S_0$ be the blowup of a surface $S_0$ at a point $p_0$. Then $S_1$ contains an exceptional curve $E_1$ (isomorphic to the projective line $\mathbb{P}^1$) that is contracted to $p_0$. Let $p_1 \in E_1$. We can blowup again at $p_1$ and choose a point $p_2$ in the new exceptional curve $E_2$. Iterating this we get a sequence of \emph{infinitely near} points $(p_0,\ldots,p_n)$.

Assuming that $S_0$ is the affine plane, then $p_0$ is defined by affine coordinates $(a_0,b_0)$ on the plane,
while for $i>0$ the point $p_i$ is defined by homogeneous coordinates $[a_i:b_i]$ on the projective line $E_i$.

The coordinates $[a_i:b_i]$ have the geometric interpretation of tangent directions (the direction of a line tangent to the branches of a curve singularity). After each blowup, the new surface $S_i$ is covered by affine plane charts, and we choose the one that contains the point $p_i$. At each step we choose coordinates such that the new exceptional curve is always the line $y=0$. More precisely, if a curve is given by $F(x,y)=0$, then its blowup is given (in the chart where the tangent direction is not $[1\!:\!0]$) by substituting $x\mapsto xy$, i.e.\ $F(xy,y)=0$; if the tangent direction is $[1\!:\!0]$, we take $F(x,xy)=0$ followed by the swap $(x,y)\mapsto (y,x)$.

In practice it is important to have this in mind in order to track the exceptional curves.

\subsection{Implementation}

The package provides the constructor
$$\texttt{LinearSystem(L,pts,m,t)},$$
where:
\begin{enumerate}
  \item[$\cdot$] $L$ is a linear system on the affine plane;
  \item[$\cdot$] $pts$ is a sequence of points in the affine plane;
  \item[$\cdot$] $m$ encodes the multiplicity sequence along the blowup chain;
  \item[$\cdot$] $t$ encodes the tangent directions chosen at each step.
\end{enumerate}

The algorithm proceeds iteratively:
\begin{enumerate}
  \item impose an ordinary multiplicity at the starting point;
  \item blowup the plane at this point
        and express the strict transform in new coordinates;
  \item divide by the exceptional factor and impose the next multiplicity (at the point given by the tangent direction);
  \item repeat until all singularities are resolved;
  \item finally blow-down to return to the original coordinates.
\end{enumerate}

Each step involves only linear algebra (evaluation matrices and nullspaces),
so the procedure is effective and scales well.

\subsection{Examples}

\noindent\textbf{Tacnode and cusp.}
A tacnode is characterized by two infinitely near double points (multiplicities \([2,2]\)) with a single tangent direction. A cusp is a double point whose strict transform becomes smooth after one blowup and is tangent to the exceptional divisor. This corresponds to multiplicities \([2,1,1]\) with two specified tangent directions, the second chosen to enforce tangency to the exceptional curve. Let us compute quartic curves with one tacnode and one cusp:
\begin{lstlisting}[language=Magma]
A<x,y>:=AffineSpace(Rationals(),2);
J:=LinearSystem(A,4);
p:=[A![0,0],A![2,3]];
m:=[[2,2],[2,1,1]];
t:=[[[1,1]],[[1,1],[1,0]]];
L:=LinearSystem(J,p,m,t);
\end{lstlisting}

\medskip
One can check the result:
\begin{lstlisting}[language=Magma]
C:=Curve(A,&+Sections(L));
[ResolutionGraph(C,q):q in p];
\end{lstlisting}

\bigskip
\noindent\textbf{A pencil of sextics.}
Now we wish to construct a pencil of plane sextic curves with nine infinitely–near double points
(multiplicity sequence \([2,\ldots,2]\) of length \(9\)) at the origin.
If we randomly fix the nine points (i.e. one point and 8 tangent directions), typically we get a double cubic.
To avoid this degeneracy, we fix only the first eight infinitely–near points and search
for the ninth one.

Let $[1,a]$ be the 8th tangent direction. Working over finite fields $\mathbb F_{p^2}$,
we check all possibilities for $a$.
For each prime \(p\) we get two values $a_1,a_2$ of the parameter \(a\) that produce a genuine
pencil. 

\begin{lstlisting}[language=Magma]
p:=59;
K:=GF(p,2);
A<x,y>:=AffineSpace(K,2);
J:=LinearSystem(A,6);
p:=A![0,0];
m:=[2,2,2,2,2,2,2,2,2]; 
for a in Set(K) diff {0} do
  t:=[[1,1],[1,2],[1,3],[1,4],[1,5],[1,6],[1,7],[1,a]];
  L:=LinearSystem(J,p,m,t);
  if Nsections(L) eq 2 then
    C:=Curve(A,&+Sections(L));
    a,ResolutionGraph(C,p);
  end if;
end for;
\end{lstlisting}

\medskip
This suggests that the desired real number $a$ is given by a quadratic extension of the rationals.
To find this extension one just needs to consider, for many primes $p$, the polynomial
\[
P(x)=(x-a_{1})(x-a_{2})
\]
and then, via CRT/RationalReconstruction in Magma, lift its coefficients to characteristic zero.
The result is $$P(x)=x^2-\frac{3645985316400}{227892834937}x+\frac{14582741040000}{227892834937}.$$

\bigskip
\noindent\textbf{Quadrifolium.}
Finally, let's compute the {\em quadrifolium}: a curve of degree 6 with a quadruple point that resolves to two different double points after one blowup. Thus it looks like the union of two tacnodes with different tangent directions. We consider this as two different singularities of type $[4, 2]$ at the same point. In order to get a nicer picture, we consider curves symmetric with respect to the coordinate axes. We also ask that the curve contains, with multiplicity 1, three additional points, one of these with tangent directions $[[1, 1]]$ (we are fixing the tangent line at the point). When we do not impose a tangent direction, the corresponding sequence of directions is the empty one: $[\ ].$ 
\begin{lstlisting}[language=Magma]
A<x,y>:=AffineSpace(Rationals(),2);
s6:=Sections(LinearSystem(A,6));
s:=[q:q in s6 | q eq Evaluate(q,[-x,y]) and q eq Evaluate(q,[x,-y])];
J:=LinearSystem(A,s);
p:=[A![0,0],A![0,0],A![1,1],A![2/10,7/10],A![7/10,2/10]];
m:=[[4,2],[4,2],[1,1],[1],[1]];
t:=[[[1,0]],[[0,1]],[[1,-1]],[],[]];
Sections(LinearSystem(J,p,m,t))[1];

\end{lstlisting}

The output is:
\begin{lstlisting}[language=Magma]
x^6+26171/9604*x^4*y^2+26171/9604*x^2*y^4-35775/4802*x^2*y^2+y^6
\end{lstlisting}

\begin{figure}[H]
  \centering
  \includegraphics[width=.5\linewidth]{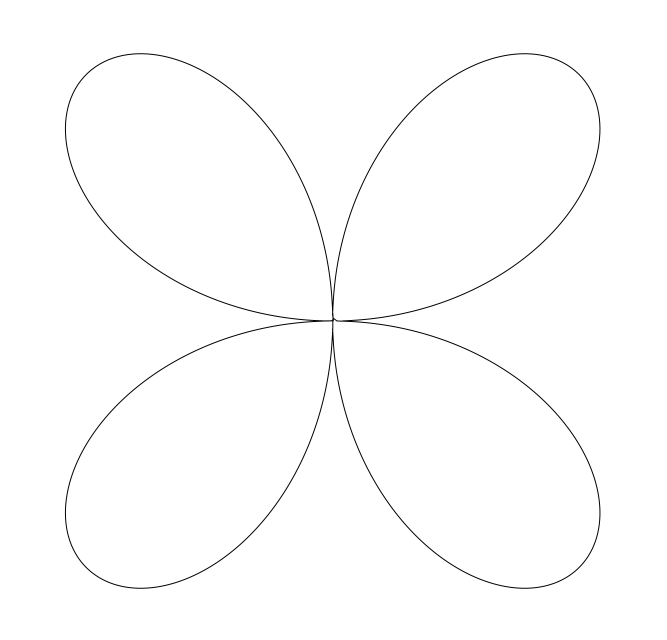}
  \caption*{Quadrifolium}
\end{figure}

\section{Examples: singular quintic surfaces}
In this section we present brief examples that are, in our view, both fun and powerful. Our aim is not novelty but simplicity: the \texttt{LinearSystem} tools make it straightforward (at least over finite fields) to write down quintic surfaces with many nodes or cusps. On the nodal side, we routinely reach $30$ nodes and also attain $31$, the sharp maximum for quintics (Beauville \cite{B}).
On the cuspidal side, we obtain examples of $\mathbb{Z}/6$–invariant surfaces with $15$ cusps plus $3$ nodes,
close to the still-unrealized target of $18$ cusps.

\subsection{$\mathbb{Z}/5$-invariant quintics with 20 nodes}

Here we construct a 4-parameter family of $\mathbb{Z}/5$-invariant quintics with 20 nodes.

We work in $\mathbb{P}^{3}$ with the $\mathbb{Z}/5$--action
\[
(x_{1}:x_{2}:x_{3}:x_{4}) \longmapsto (x_{1}:r^{2}x_{2}:r x_{3}:r x_{4}), \qquad r^{5}=1.
\]
The corresponding invariant quintics form a $13$-dimensional subspace spanned by the monomials listed in the code below. Over the function field $\mathbb F=\mathbb{Q}(a,b,c,d)$ we build the linear system of invariant quintics and impose double points at four points,
\[
p_{1}=(1\!:\!1\!:\!1\!:\!1),\quad
p_{2}=(3\!:\!3\!:\!2\!:\!1),\quad
p_{3}=(a\!:\!a\!:\!b\!:\!1),\quad
p_{4}=(c\!:\!c\!:\!d\!:\!1).
\]
Computer experiments suggested that points with $x_{1}=x_{2}$ impose one fewer condition on the invariant subspace,
so we deliberately choose representatives with the first two coordinates equal.
For general parameters the  $\mathbb{Z}/5$--orbits have size $5$, yielding a surface with $4\times 5=20$ ordinary double points.
The Magma code below constructs this system and returns one polynomial $$F=F_{a,b,c,d}(x_1,x_2,x_3,x_4).$$
We clear denominators at the end and save the polynomial in a file.

\medskip

\begin{lstlisting}[language=Magma]
F<a,b,c,d>:=FunctionField(Rationals(),4);
P3<x1,x2,x3,x4>:=ProjectiveSpace(F,3);
s:=[x1^5,x2^5,x1^2*x2^2*x3,x1*x2*x3^3,x3^5,x1^2*x2^2*x4,x1*x2*x3^2*x4,
    x3^4*x4,x1*x2*x3*x4^2,x3^3*x4^2,x1*x2*x4^3,x3^2*x4^3,x3*x4^4,x4^5];
L:=LinearSystem(P3,s);
L:=LinearSystem(L,[P3![1,1,1,1],P3![3,3,2,1]],[2,2]);
L:=LinearSystem(L,[P3![a,a,b,1],P3![c,c,d,1]],[2,2]);
F:=Sections(L)[1];
lcm:=LCM({Denominator(q):q in Coefficients(F)});
F:=lcm*F;
\end{lstlisting}

\medskip
This took only 0.2 seconds in our computer!

\subsection{Quintics with 30 or 31 nodes}

Working over a quadratic extension $\mathbb{F}_{p^{2}}$, we specialize at random the parameters
$(a,b,c,d)$ in the polynomial $F$ computed above, which defines a $\mathbb Z/5$--invariant
four--parameter family of $20$–nodal quintics. For each specialization we form the surface
$X\subset\mathbb{P}^{3}$ and compute its singular subscheme $S={\rm Sing}\, (X)$. The search is
fully automatic: we simply test whether
\[
\deg(S)=\deg\bigl(\mathrm{ReducedSubscheme}(S)\bigr)\in\{30,31\}.
\]
In practice, this procedure quickly produces examples
with $30$ nodes and even with $31$ nodes, the sharp maximum for quintics (Beauville):

\medskip

\begin{lstlisting}[language=Magma]
K:=GF(101);
P3<x1,x2,x3,x4>:=ProjectiveSpace(K,3);
F:=x1^5+x2^5+76*x1^2*x2^2*x3+54*x1*x2*x3^3+65*x3^5+90*x1^2*
   x2^2*x4+93*x1*x2*x3^2*x4+29*x3^4*x4+37*x1*x2*x3*x4^2+53*
   x3^3*x4^2+85*x1*x2*x4^3+20*x3^2*x4^3+10*x3*x4^4+93*x4^5;
G:=x1^5+x2^5+48*x1^2*x2^2*x3+62*x1*x2*x3^3+97*x3^5+5*x1^2*
   x2^2*x4+90*x1*x2*x3^2*x4+12*x3^4*x4+80*x1*x2*x3*x4^2+99*
   x3^3*x4^2+61*x1*x2*x4^3+36*x3^2*x4^3+18*x3*x4^4+97*x4^5;
X:=Scheme(P3,F); Y:=Scheme(P3,G);
SX:=ReducedSubscheme(SingularSubscheme(X));
SY:=ReducedSubscheme(SingularSubscheme(Y));
Degree(SX) eq 30;
Degree(SY) eq 31;
\end{lstlisting}

\medskip
With a little additional work, these finite–field examples can be lifted to characteristic
zero, but we do not pursue this here.

\subsection{$\mathbb{Z}/6$-invariant quintics with 15 nodes}

Here we construct a 6-parameter family of $\mathbb{Z}/6$-invariant quintics with 15 nodes.

The action is given by
\[
\begin{aligned}
(x_1:x_2:x_3:x_4) &\mapsto (x_1:x_2:-x_3:x_4)\\
(x_1:x_2:x_3:x_4) &\mapsto (r^2 x_1: r x_2: x_3: x_4),
\end{aligned}
\]
with \(r^3=1\).
The corresponding invariant quintics form a $11$-dimensional subspace spanned by the monomials listed in the code below.
We impose one \(\mathbb{Z}/2\)-fixed ordinary double point at \([1:1:0:1]\), which leaves space to imposing two further
general double points. This yields a 6-parameter family, given by a single polynomial
$$F=F_{a,b,c,d,e,f}(x_1,x_2,x_3,x_4),$$ whose general member has exactly $15$ nodes.

\medskip
\begin{lstlisting}[language=Magma]
K:=Rationals();
F<a,b,c,d,e,f>:=FunctionField(K,6);
R<x1,x2,x3,x4>:=PolynomialRing(F,4,"grevlex");
P3:=ProjectiveSpace(R);
s5:=[x4^5,x3^2*x4^3,x1*x2*x4^3,x2^3*x4^2,x1^3*x4^2,x3^4*x4,x1*x2*x3^2*x4,
     x1^2*x2^2*x4,x2^3*x3^2,x1^3*x3^2,x1*x2^4,x1^4*x2];
L:=LinearSystem(P3,s5);
L:=LinearSystem(L,P3![1,1,0,1],2);
L:=LinearSystem(L,[P3![a,b,c,1],P3![d,e,f,1]],[2,2]);
F:=Sections(L)[1];
lcm:=LCM([Denominator(q):q in Coefficients(F)]);
F:=lcm*F;
\end{lstlisting}

\medskip
This took only 0.2 seconds in our computer!

\subsection{Quintics with many cusps}
Analogously to Section~5.2, we performed a random search within the previously computed $\mathbb{Z}/6$–invariant family of quintic surfaces, over finite fields. This gave examples with $15$ cusps, as well as examples with $15$ cusps plus $3$ nodes.

\medskip
\begin{lstlisting}[language=Magma]
K:=GF(103);
P3<x1,x2,x3,x4>:=ProjectiveSpace(K,3);
F:=x1^4*x2+30*x1*x2^4+22*x1^3*x3^2+29*x2^3*x3^2+85*x1^2*x2^2*x4+
   25*x1*x2*x3^2*x4+56*x3^4*x4+15*x1^3*x4^2+89*x2^3*x4^2+
   60*x1*x2*x4^3+22*x3^2*x4^3+29*x4^5;
G:=x1^4*x2+42*x1*x2^4+73*x1^3*x3^2+60*x1^2*x2^2*x4+
   9*x1*x2*x3^2*x4+93*x3^4*x4+15*x1^3*x4^2+77*x2^3*x4^2+
   98*x1*x2*x4^3+39*x3^2*x4^3+16*x4^5;
X:=Surface(P3,F); Y:=Surface(P3,G);
ptsX:=SingularPoints(X);
ptsY:=SingularPoints(Y);
#ptsX eq 15, #ptsY eq 18;
for q in ptsX do IsSimpleSurfaceSingularity(X!q);end for;
for q in ptsY do IsSimpleSurfaceSingularity(Y!q);end for;
\end{lstlisting}

\vspace{1cm}

\noindent Carlos Rito
\vspace{0.1cm}
\\ Centro de Matem\'atica, Universidade do Minho - Polo CMAT-UTAD
\vspace{0.1cm}
\\ Universidade de Tr\'as-os-Montes e Alto Douro, UTAD
\\ Quinta de Prados
\\ 5000-801 Vila Real, Portugal
\vspace{0.1cm}
\\ www.utad.pt, {\tt crito@utad.pt}

\end{document}